\documentclass[11pt]{article}
\usepackage{}
\usepackage{dsfont}
\usepackage{amsmath}
\usepackage{mathrsfs}
\usepackage{amsmath,amssymb}
\usepackage{amsfonts}
\usepackage{hyperref}
\usepackage{amsthm}
\usepackage{graphicx}
\renewcommand{\qed}{\hfill\small{$\square$}\normalsize}

\theoremstyle{definition}
\newtheorem{lemma}{Lemma}[section]
\newtheorem{definition}[lemma]{Definition}

\newtheorem{theorem}[lemma]{Theorem}
\newtheorem{corollary}[lemma]{Corollary}
\newtheorem{example}{Example}
\newtheorem{remark}{Remark}

\numberwithin{equation}{section}
\renewcommand{\proof}{\textbf{Proof. }}
\renewcommand{\qed}{\hfill\small{$\square$}\normalsize}

\DeclareFixedFont{\Acknowledgment}{OT1}{cmr}{bx}{n}{18pt}
\begin{document}

\title{\bf Discrete Quasi-Einstein Metrics and Combinatorial Curvature Flows in 3-Dimension}
\author{Huabin Ge, Xu Xu}
\maketitle

\abstract{}
We define Discrete Quasi-Einstein metrics (DQE-metrics) as critical points of discrete total curvature functional on triangulated 3-manifolds. We study DQE-metrics by introducing combinatorial curvature flows. We prove that these flows produce solutions which converge to discrete quasi-Einstein metrics when the initial energy is small enough. The proof relies on a careful analysis of discrete dual-Laplacians which we interpret as the Jacobian matrix of the curvature map. As a consequence, combinatorial curvature flow provides an algorithm to compute discrete sphere packing metrics with prescribed curvatures.

\section{Introduction}
In \cite{CR}, Cooper and Rivin introduced a type of combinatorial scalar curvature which is defined as angle defect of solid angles. Their curvature is a combinatorial analogy of the scalar curvature in the infinitesimal sense and is closely related to Regge's calculus. They used total curvature functional to study the combinatorial analogue of conformal deformation of sphere packing metrics. They found that the Hessian of total curvature functional is crucial.

In \cite{G1}, Professor Glickenstein introduced a combinatorial version of Yamabe flow. The flow is a system of equations connecting CR-curvature (Cooper-Rivin's curvature) to sphere packing metrics. The evolution of CR-curvature was then derived to satisfy a combinatorial heat equation. The discrete Laplacian was defined so as to describe the combinatorial heat equation. Moreover, \cite{G1} studied the dual structure in Euclidean triangulations coming from sphere packing, and expressed the discrete Laplacian by the geometric elements in the dual structure.

In this paper, we define DQE-metric (discrete quasi-Einstein metric) as the critical point of total curvature functional and give some analytical conditions for the existence of DQE-metrics. We study DQE-metric by combinatorial curvature flow methods. This paper is mainly concerned with two types of discrete energies and their negative gradient flows. One is the total curvature functional which determines a combinatorial CR-curvature flow of second order. The other is the discrete quadratic energy, a 3-dimensional analogue of discrete Calabi energy introduced in \cite{Ge}, which determines a combinatorial CR-curvature flow of fourth order. These flows exhibit good existence and convergence properties, that is, solutions of these flows converge to discrete quasi-Einstein metrics when the initial normalized quadratic energy is small enough. Thus these flows are useful tools for the study of manifolds with sphere packing metrics. These flows can be modified as user's will so as the modified flows tend to find user prescribed target curvatures automatically. Thus we can design algorithm to compute discrete sphere packing metrics with user prescribed curvatures. It's interesting that much of our process can be formally generalized to higher dimensions and to other definitions of curvatures.

\section{Conformal structure and DQE-metric}
\subsection{Sphere packing metric and edge-tangential-sphere}
Consider a compact manifold $M$ of dimension 3 with a triangulation $\mathcal{T}$ on $M$. The triangulation is written as $\mathcal{T}=\{V,E,F,T\}$, where the symbols $V,E,F,T$ represent the set of vertices, edges, faces and tetrahedrons respectively. Throughout this paper, all functions defined on vertices are regarded as a column vector and $N=V^{\#}$, the number of vertices. Furthermore, all vertices are marked by $v_{1},\cdots,v_{N}$ one by one. We often write $i$ instead of $v_i$.

Sphere packing metric is a map $r:V\rightarrow (0,+\infty)$ such that the length between vertices $i$ and $j$ is $l_{ij}=r_{i}+r_{j}$ for each edge $\{i,j\}\in E$,
and the lengths $l_{ij},l_{ik},l_{il},l_{jk},l_{jl},l_{kl}$ can be realized as an Euclidean or hyperbolic tetrahedron for each tetrahedron $\{i,j,k,l\}\in T$. In \cite{CR}, Cooper-Rivin called tetrahedrons generated by this way ``conformal". Mathematicians had found that a tetrahedron is conformal iff there exists a unique sphere tangent to all of the edges of the tetrahedron. Moreover, the point of tangency with the edge $\{i,j\}$ is at distance $r_i$ from $v_i$.

We may think of sphere packing metrics as points in $\mathds{R}^N_{>0}$, $N$ times Cartesian product of $(0,\infty)$. However, not all points in $\mathds{R}^N_{>0}$ represent sphere packing metrics, for we need extra nondegenerate conditions besides positive condition. A conformal tetrahedron $\{i,j,k,l\}\in T$ generated by four positive radii $r_{i},r_{j},r_{k},r_{l}$ can be realized as an Euclidean tetrahedron of same edge lengths if and only if $$Q_{ijkl}=\left(\frac{1}{r_{i}}+\frac{1}{r_{j}}+\frac{1}{r_{k}}+\frac{1}{r_{l}}\right)^2-
2\left(\frac{1}{r_{i}^2}+\frac{1}{r_{j}^2}+\frac{1}{r_{k}^2}+\frac{1}{r_{l}^2}\right)>0.$$
thus, the space of admissible sphere packing metrics is $$\mathfrak{M}_{r}=\left\{\;r\in\mathds{R}^N_{>0}\;\big|\;Q_{ijkl}>0, \;\forall \{i,j,k,l\}\in T\;\right\}.$$
In \cite{CR}, Cooper and Rivin showed that $\mathfrak{M}_{r}$ is an open set in $\mathds{R}^N_{>0}$, simplex connected, but not convex.
\remark For a $n$-simplex, conformal condition, $i.e.$, the existence of sphere packing metric, is equivalent to the existence of a $n$ dimensional sphere tangent to each of the edges.

\subsection{CR-curvature and total curvature functional}
Given a Euclidean tetrahedron $\{i,j,k,l\}\in T$, the solid angle of vertex $i$ is denoted by $\alpha_{ijkl}$. Cooper and Rivin's combinatorial scalar curvature is defined as angle defect of solid angles. Concretely, the curvature $K_{i}$ at a vertex $i$ is $$K_{i}\doteqdot 4\pi-\sum_{\{i,j,k,l\}\in T}\alpha_{ijkl},$$
here, the sum is taken over all tetrahedrons having $i$ as one of its vertex. We write this curvature as ``CR-curvature" for short. The curvature $K_{i}$ is a combinatorial analogue of scalar curvature in smooth cases, for it measures the difference between, on the one hand, the growth rate of a small ball centered at vertex $\{i\}$ in $M$ and, on the other hand, the growth rate of the volume of an Euclidean ball of the same radius.

Cooper and Rivin first defined this type of combinatorial curvature in \cite{CR}, they also studied the sphere packing metric structure and proved that the metric cannot be deformed while keeping the solid angles fixed. Their central conclusion were derived from a careful study of the following total curvature functional $$S=\sum_{i=1}^N K_{i}r_{i}.$$ They showed that the functional is weakly concave and is strongly concave up to scaling. In fact, by the famous Schl$\ddot{a}$fli formula, one can get $$\nabla_{r}S=K,$$ $$Hess_{r}S=\frac{\partial(K_{1},\cdots,K_{N})}{\partial(r_{1},\cdots,r_{N})}.$$ From \cite{CR}, \cite{Ri}, \cite{G2} or Lemma \ref{semi-positive} proved in this paper, the curvature map $K=K(r)=\nabla_{r}S$ is locally an embedding when restricted to the space transverse to the scaling direction $r$. Equivalently, the sphere packing metric can be determined by the combinatorial scalar curvature locally up to scaling.

Glickenstein introduced the so called combinatorial Yamabe flow $$\frac{dr_{i}}{dt}=-K_{i}r_{i}$$ to study combinatorial Yamabe problem in \cite{G1}. This geometric flow connects CR-curvatures and sphere packing metrics. Glickenstein showed that the evolution of CR-curvature satisfies a heat equation driven by discrete dual-Laplacian.

\cite {CL1}, \cite{CR} and \cite{G1} provided great perspectives and deep inspirations for the study of combinatorial manifolds with circle packing metrics. The combinatorial structure of the triangulation consists of finite vertices, edges, faces, tetrahedrons and their adjacent relations. We can always represent all the finite relations by incidence matrix. Methodologically, it's more convenient to use the language of matrix when dealing problems arisen in combinatorial manifolds. Under this belief, we will give a new interpretation of discrete dual-Laplacian as the Jacobian matrix of the curvature map.

\section{Discrete quasi-Einstein metric}
Finding good metrics is always a central topic in geometry. For a smooth manifold, Einstein metric is a good candidate for privileged metric on the manifold. On an Einstein manifold, the Riemannian metric $g$ is proportionate to the Ricci curvature, that is $Ric=\lambda g$.

For a combinatorial manifold $(M,\mathcal{T})$, consider sphere packing metric $r$ as an analogue of Riemannian metric $g$, while CR-curvature $K$ as Ricci curvature $Ric$, we write $K=\lambda r$ formally as the analogue of Einstein condition $Ric=\lambda g$.
\definition Given $(M^3,\mathcal{T})$, where $M^3$ is a compact manifold, $\mathcal{T}$ is a fixed triangulation on $M^3$. A sphere packing metric $r$ with $K=\lambda r$ is called a \textbf{DQE-metric} (discrete quasi-Einstein metric). A DQE-metric $K=\lambda r$ is called ``\textbf{flat}", ``\textbf{positive}" or ``\textbf{negative}" when the constant $\lambda$ satisfies $\lambda=0$, $\lambda>0$ or $\lambda<0$ respectively.\\[0.1pt]

Notice that, when $r$ is a DQE-metric, \emph{i.e.} $K=\lambda r$, then $\lambda=\frac{S}{\|r\|^2}$. Besides their formal similarity to smooth Einstein metrics, DQE-metrics themselves have great importance. They are the critical points of the normalized total curvature functional.

\theorem Consider the $\tau$-normalized total curvature functional $S_\tau=\frac{S}{\|r\|^\tau}$.
\begin{description}
  \item[(1)] If $\tau=1$, the critical points of $S_\tau$ are DQE-metrics, vice versa.
  \item[(2)] If $\tau\neq1$, the critical points of $S_\tau$ are flat DQE-metrics, vice versa.
\end{description}
\proof
\begin{equation*}
\nabla_{r}S_{\tau}=\frac{1}{\|r\|^{\tau}}(\nabla_{l}S-\frac{\tau S}{\|r\|^{2}}r)=\frac{1}{\|r\|^{\tau}}(K-\frac{\tau S}{\|r\|^{2}}r).
\end{equation*}
$r$ is the critical point of $S_\tau$ if and only if $K=\frac{\tau S}{\|r\|^{2}}r$. Thus the conclusions in the theorem follow from $S=r^TK=\tau S$. \qed

We want to know whether DQE-metric exists on a fixed $(M,\mathcal{T})$. Inspired by Thurston's work in \cite{T1}, we believe that there are combinatorial and topological obstructions for the existence of DQE-metrics.

Provided combinatorial curvature flow methods (\cite{CL1}, \cite{G1}), we study DQE-metrics by introducing combinatorial curvature flows in the following section. These flows are negative gradient flows of some discrete energies. We are mainly concerned about two types of discrete energies and their negative gradient flows. One is the total curvature functional $S=r^TK$, which determines a combinatorial CR-curvature flow of second order. The other is the discrete quadratic energy $\mathcal{C}=\|K\|^2$, which determines a combinatorial CR-curvature flow of fourth order. The latter seems more powerful than the former, so we introduce the fourth order flow first.

\section{Combinatorial CR-curvature flow of fourth order} \label{4th}
\subsection{Definition}
Consider the discrete quadratic energy functional:
\begin{equation}
\mathcal{C}(r)=\|K\|^2=\sum_{i=1}^N K_i^2
\end{equation}
Choose coordinate transformations $u_i=\ln r_i$, which is a homeomorphism from space $\mathfrak{M}_r$ to space $\mathfrak{M}_u$. We do not distinguish the quadratic energy functional as a function of $r$ or as a function of $u$.People can distinguish them from the context without confusion. Remember, all variables $r$, $u$ and $K$ defined on vertices are considered as $N$-dimensional column vector.

Under $u$-coordinate, the Jacobian matrix of CR-curvature map $K=K(u)$ is
\begin{displaymath}
L=(L_{ij})_{N\times N}=\frac{\partial(K_{1},\cdots,K_{N})}{\partial(u_{1},\cdots,u_{N})}=
\left(
\begin{array}{ccccc}
 {\frac{\partial K_1}{\partial u_1}}& \cdot & \cdot & \cdot &  {\frac{\partial K_1}{\partial u_N}} \\
 \cdot & \cdot & \cdot & \cdot & \cdot \\
 \cdot & \cdot & \cdot & \cdot & \cdot \\
 \cdot & \cdot & \cdot & \cdot & \cdot \\
 {\frac{\partial K_N}{\partial u_1}}& \cdot & \cdot & \cdot &  {\frac{\partial K_N}{\partial u_N}}
\end{array}
\right),
\end{displaymath}
then we have
\begin{displaymath}
\nabla_u \mathcal{C}=
\left(
\begin{array}{c}
 {\frac{\partial \mathcal{C}}{\partial u_1}} \\
 \cdot\\
 \cdot\\
 \cdot\\
 {\frac{\partial \mathcal{C}}{\partial u_N}}
\end{array}
\right)
=2
\left(
\begin{array}{ccccc}
 {\frac{\partial K_1}{\partial u_1}}& \cdot & \cdot & \cdot &  {\frac{\partial K_N}{\partial u_1}} \\
 \cdot & \cdot & \cdot & \cdot & \cdot \\
 \cdot & \cdot & \cdot & \cdot & \cdot \\
 \cdot & \cdot & \cdot & \cdot & \cdot \\
 {\frac{\partial K_1}{\partial u_N}}& \cdot & \cdot & \cdot &  {\frac{\partial K_N}{\partial u_N}}
\end{array}
\right)
\left(
\begin{array}{c}
 K_1 \\
 \cdot\\
 \cdot\\
 \cdot\\
 K_N
\end{array}
\right)=2L^TK.
\end{displaymath}

\definition Combinatorial CR-curvature flow of fourth order is defined to be the following gradient flow
\begin{equation} \label{Combina-Calabi-F}
\dot{u}(t)=-\frac{1}{2}\nabla_{u} \mathcal{C}=-L^TK,
\end{equation}
where ``$\cdot$" means time derivative, and the derivative of column vector means derivative of its components.

The CR-curvature evolves according to $\dot{K}=-LL^TK$. $LL^T$ is a discrete fourth order differential operator which acts on function (defined on vertices) by matrix multiplication. This is why flow (\ref{Combina-Calabi-F}) is of fourth order in our denotion.

We call $L$ discrete dual-Laplacian and write ``DDL" for short. It is somewhat abrupt to give $L$ this name, however, $L$ is exactly a type of discrete Laplace operator whose definition relies on the dual structure of sphere packing metric and conformal tetrahedrons. Next we will look at $L$ intensively and interpret how the name ``DDL" came about.

\subsection{Discrete dual-Laplacian}
Regard the 1-skeleton of the triangulation as a 3 dimensional graph. In spectral graph theory, discrete Laplacians are always in the following form (see \cite{CHU})
$$\Delta f_{i}=\sum_{j\thicksim i}B_{ij}(f_{j}-f{i}),$$ where $B_{ij}$ are weights along each edge $\{i,j\}$. ``$j\thicksim i$" means two vertices $i$ and $j$ are adjacent, and the sum is taken over every edge having $i$ as one of its vertex. Weights function $B_{\bullet\bullet}$ gives each edge $j\thicksim i$ a weight $B_{ij}$, different weights function give different types of discrete Laplacians. Conformal tetrahedron has a strong dual structure, which determines a special weights function and the corresponding DDL.

For adjacent vertices $i\thicksim j$, the $(i,j)$-entry of $L$ is
$$\frac{\partial K_{i}}{\partial u_{j}}=-\frac{\partial \left(\sum_{\{i,j,k,l\}\in T}\alpha_{ijkl}\right)}{\partial u_{j}}=-\sum_{\{i,j,k,l\}\in T}\frac{\partial\alpha_{ijkl}}{\partial u_{j}}=-\sum_{\{i,j,k,l\}\in T}\frac{\partial\alpha_{ijkl}}{\partial r_{j}}r_{j}.$$
Let
$$B_{ij}=-\frac{\partial K_{i}}{\partial u_{j}}=\sum_{\{i,j,k,l\}\in T}\frac{\partial\alpha_{ijkl}}{\partial r_{j}}r_{j}.$$
\lemma \label{L-B-relation}
For any $1\leq i,j\leq N$,
\begin{gather}
L_{ij}=
\begin{cases}
\sum\limits_{k \sim i}B_{ik} \,, & \text{$ j=i $} \\
-B_{ij} \,,& \text{$ j\sim i $}\;\\
0 \,,& \text{$ else $}
\end{cases}
\end{gather}\\[5pt]
\proof Solid angles are homogeneous functions of degree zero of sphere packing metrics. Use Euler's homogeneous function theorem, we get
\begin{equation}\label{0-homogeneous}
r_{i}\frac{\partial\alpha_{ijkl}}{\partial r_{i}}+r_{j}\frac{\partial\alpha_{ijkl}}{\partial r_{j}}+r_{k}\frac{\partial\alpha_{ijkl}}{\partial r_{k}}+r_{l}\frac{\partial\alpha_{ijkl}}{\partial r_{l}}=0,
\end{equation}
which implies $L_{ii}=B_{ii}=\sum\limits_{j\thicksim i}B_{ij}$. \qed
\remark Glickenstein had already found above formula \ref{0-homogeneous} using Schl$\ddot{a}$fli formula in \cite{G1}.

\corollary $L(1,\cdots,1)^T=0.$\\[2pt]

Given a conformal 3-simplex $\tau=\{i,j,k,l\}$.
\begin{figure}
\centering
 \unitlength 1mm 
\linethickness{0.4pt}
\ifx\plotpoint\undefined\newsavebox{\plotpoint}\fi 
\begin{picture}(92,55)(0,0)
\thicklines
\multiput(58.5,64)(-.0337398374,-.0788617886){615}{\line(0,-1){.0788617886}}
\multiput(37.75,15.5)(-.0336787565,.0446891192){386}{\line(0,1){.0446891192}}
\multiput(24.75,32.75)(.03640776699,.03371089536){927}{\line(1,0){.03640776699}}
\multiput(58.5,64)(.0336842105,-.0610526316){475}{\line(0,-1){.0610526316}}
\multiput(74.5,35)(-.0631487889,-.0337370242){578}{\line(-1,0){.0631487889}}
\multiput(60.75,28)(.03125,.03125){8}{\line(0,1){.03125}}
\multiput(24.68,32.43)(.995,.055){51}{{\rule{.8pt}{.8pt}}}
\multiput(61,27.75)(-.03370787,.12921348){89}{\line(0,1){.12921348}}
\put(61.5,24.25){D}
\put(76.75,35){j}
\put(36.75,12){i}
\put(60.25,65){l}
\put(20,32){k}
\multiput(51.18,36.68)(0,-.90625){9}{{\rule{.8pt}{.8pt}}}
\multiput(57.93,38.93)(-.8125,-.28125){9}{{\rule{.8pt}{.8pt}}}
\multiput(51.18,28.93)(.95,-.125){11}{{\rule{.8pt}{.8pt}}}
\put(58,40){F}
\put(48.25,36.5){C}
\put(47.75,27){E}
\end{picture}

 \caption{Dual structure in conformal tetrahedron}\label{dual-structure}
\end{figure}
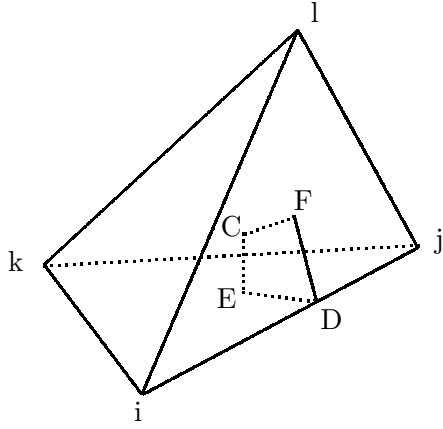
Let $C$ denote the center of its edge-tangential-sphere. $E$ and $F$ are foots of perpendicular from $C$ to face $\{i,j,k\}$ and face $\{i,j,l\}$ respectively. The 2-plane spanned by $C,E,F$ meets edge $\{i,j\}$ perpendicularly at point $D$. Then $E$ ($F$) is the center of circle inscribed in the Euclidean triangle $\{i,j,k\}$ ($\{i,j,l\}$). In a single conformal 3-simplex $\tau$, the dual $\star\{i,j\}_\tau$ to the edge $\{i,j\}$ (see \cite{G1}) is defined to be the region surrounded by four points $C,E,D,F$. The measurement of $\star\{i,j\}_\tau$ is the directed area of this region, which is denoted by $A_{ijkl}$ (See Figure \ref{dual-structure}). The area is negative when the center of edge-tangential-sphere is at the outside of the conformal tetrahedron, and positive when at inside. In the whole triangulation $(M,\mathcal{T})$ with fixed sphere packing metric, the dual $\star\{i,j\}$ is the sum of all such $\star\{i,j\}_\tau$, where $\tau$ runs over all conformal tetrahedrons containing $\{i,j\}$ as one edge. Denote $l_{ij}^\ast$ as the measurement of $\star\{i,j\}$, then
$$l_{ij}^\ast=\sum_{\{i,j,k,l\}\in T}A_{ijkl}.$$
Glickenstein proved (\cite{G1})
$$\frac{\partial\alpha_{ijkl}}{\partial r_{j}}r_{j}r_{i}=\frac{A_{ijkl}}{l_{ij}}.$$
Therefore $$B_{ij}=\sum_{\{i,j,k,l\}\in T}\frac{\partial\alpha_{ijkl}}{\partial r_{j}}r_{j}=\frac{1}{r_{i}}\frac{l_{ij}^\ast}{l_{ij}}.$$
It implies that $$r_{i}B_{ij}=\frac{l_{ij}^\ast}{l_{ij}}=r_{j}B_{ji}$$ is symmetric in the two indices.

Glickenstein defined a discrete Laplace operator $``\Delta"$ in \cite{G1} and showed that
$$\Delta f_{i}=\frac{1}{r_{i}}\sum_{j\thicksim i}\frac{l_{ij}^\ast}{l_{ij}}(f_{j}-f_{i}).$$
Using Lemma \ref{L-B-relation} we obtain
$$\Delta f_{i}=\sum_{j\thicksim i}B_{ij}(f_{j}-f_{i})=-\left(L_{ii}f_{i}+\sum_{j=1,j\neq i}^NL_{ij}f_{j}\right)=(-Lf)_{i}.$$
Regard the discrete Laplace operator $\Delta$ as a matrix (it acts on function $f$ by matrix multiplication). Above formula is just the components form of the following $$-\Delta=L=\frac{\partial(K_{1},\cdots,K_{N})}{\partial(u_{1},\cdots,u_{N})}.$$
Thus the discrete Laplace operator $``-\Delta"$ is interpreted as the Jacobian of CR-curvature map.

Denote $R=diag(r_{1},...,r_{N})$, $\widetilde{L}=RL$. $\widetilde{\Delta}=R\Delta=-\widetilde{L}$ is called weighted DDL. Let $\Lambda=\frac{\partial(K_{1},\cdots,K_{N})}{\partial(r_{1},\cdots,r_{N})}$, then $L=\Lambda R$ and $\widetilde{L}=RL=R\Lambda R$.

\lemma \label{semi-positive}
The matrix $\widetilde{L}=RL$ is positive semi-definite, having rank $N-1$. The kernel of $\widetilde{L}$ is the span of the vector $(1,\cdots,1)^T$. Moreover, the matrix equation $L^Tx=0$ have a unique nonzero solution $x=tr$ up to scaling, where $t$ can be any constant.\\
\proof It's sufficient to prove that $\Lambda$ is positive semi-definite. For any conformal tetrahedron $\{i,j,k,l\}\in T$, it had been shown in the appendix of \cite{G2} (or see \cite{CR} and \cite{Ri}) that the following $4\times 4$ matrix
\begin{displaymath}
\Lambda_{ijkl}=-\frac{\partial(\alpha_{i},\alpha_{j},\alpha_{k},\alpha_{l})}{\partial(r_{i},r_{j},r_{k},r_{l})}=-
\left(
\begin{array}{cccc}
 {\frac{\partial \alpha_i}{\partial r_i}}& \cdot & \cdot &  {\frac{\partial \alpha_i}{\partial r_l}} \\
 \cdot & \cdot & \cdot & \cdot \\
 \cdot & \cdot & \cdot & \cdot \\
 \cdot & \cdot & \cdot & \cdot \\
 {\frac{\partial \alpha_l}{\partial r_i}}& \cdot & \cdot &  {\frac{\partial \alpha_l}{\partial r_l}}
\end{array}
\right),
\end{displaymath}
is positive semi-definite, having rank 3 and the nullspace is the span of the vector $(r_{i},r_{j},r_{k},r_{l})^T$.

We want to extend the matrix $\Lambda_{ijkl}$ to a $N\times N$ matrix. Remember that all vertices are ordered and marked by $v_{1},...,v_{N}$. We suppose $v_i$, $v_j$, $v_k$, $v_l$ arise at $\hat{i}$, $\hat{j}$, $\hat{k}$, $\hat{l}$ position respectively in the ordered sequence $v_{1},...,v_{N}$. Then we get a $N\times N$ matrix by putting $(\Lambda_{ijkl})_{st}=-\frac{\partial \alpha_s}{\partial r_t}$ at the $(\hat{s},\hat{t})$-entry position for any $s,t\in \{i,j,k,l\}$, and putting $0$ at other entries. Without confusion, we may still write the extended $N\times N$ matrix as $\Lambda_{ijkl}$,
then we have $$\Lambda=\frac{\partial(K_{1},\cdots,K_{N})}{\partial(r_{1},\cdots,r_{N})}=\sum_{\{i,j,k,l\}\in T}\Lambda_{ijkl}.$$
Because each component in the sum is positive semi-definite matrixes, $\Lambda$ is positive semi-definite .

If $\Lambda x=0$, then $\Lambda_{ijkl}x=0$, \emph{i.e.} $x\in Ker(\Lambda_{ijkl})$ for any $\{i,j,k,l\}\in T$.
Hence there is a constant $t_{ijkl}$ \emph{s.t.} $(x_{i},x_{j},x_{k},x_{l})=t_{ijkl}(r_{i},r_{j},r_{k},r_{l})$. Because the manifold $M$ is connected, all $t_{ijkl}$ must be equal. Thus $x=tr$, which implies $Ker(\Lambda)=\{tr|t\in \mathds{R}\}$ and $rank(\Lambda)=N-1$.

Since $\Lambda r=0$ and $\Lambda$ is symmetry, then $r^T\Lambda=r^TLR^{-1}=0$, hence $r^TL=0$. So $Lx=0$ has a unique solution $x=tr$ up to scaling.\qed

Due to Lemma \ref{semi-positive}, $\widetilde{L}$ has full rank in the space transverse to scaling. Also, matrix $\Lambda$ is positive definite in the space transverse to $r=(r_{1},\cdots,r_{N})^T$ direction, or equivalently, $\Lambda$ is positive definite up to scaling.
\corollary
Along the combinatorial CR-curvature flow (\ref{Combina-Calabi-F}), the quadratic energy functional is descending. Furthermore, both $\sum_{i=1}^Nu_{i}$ and $\prod_{i=1}^Nr_{i}$ are constants.
\remark By Schl$\ddot{a}$fli formula, the differential 1-form $\omega=\sum_{i=1}^NK_{i}dr_{i}=dS$ is exact. Moreover, because $\mathfrak{M}_{r}$ is simply connected, we have $$S=\int_a^r\sum_{i=1}^NK_idr_i+C,$$ where $a$ is an arbitrarily selected point form $\mathfrak{M}_{r}$.
\remark \label{transverse} For any positive constant $C>0$, $\prod_{i=1}^Nr_{i}=C$ determines a hypersurface in $\mathfrak{M}_{r}\subset\mathds{R}^N$. At each point $r$ of this hypersurface, the direction toward which $r$ scales is transverse to the hypersurface $\prod_{i=1}^Nr_{i}=C$. Thus $\widetilde{L}$ has full rank, and $\Lambda$ is positive definite along the combinatorial CR-curvature flow \ref{Combina-Calabi-F}.

\subsection{Convergence to DQE-metrics}
We had shown that DQE-metrics are critical points of $\tau$-normalized total CR-curvature functionals. Notice that, DQE-metrics are critical points of the discrete quadratic energy $\mathcal{C}$ too. Due to Lemma \ref{semi-positive}, $\nabla_{u} \mathcal{C}=2L^TK=0$ has only DQEM solutions, that is $K=\lambda r$. Next we will give an equivalent analytical condition for the existence of DQE-metrics by combinatorial CR-curvature flow methods.

\theorem If the solution to combinatorial CR-curvature flow \ref{Combina-Calabi-F} exists for all time and converges to a non-degenerate sphere packing metric $r(+\infty)$, then DQE-metric exists. Moreover, $r(+\infty)$ is indeed one.\\
\proof Along combinatorial CR-curvature flow $\dot{u}=-L^TK$, we have $$\dot{L}=-LL^TK$$ and $$\dot{\mathcal{C}}=-2K^TLL^TK\leq 0.$$ The existence of $\dot{\mathcal{C}}(+\infty)$ follows from the existence of $r(+\infty)$. This implies $\dot{\mathcal{C}}(+\infty)=0$ due to the fact that $\mathcal{C}$ is descending along flow $\dot{u}=-L^TK$. Hence $L^T(+\infty)K(+\infty)=0$, using Lemma \ref{semi-positive} we get $K(+\infty)=\lambda r(+\infty)$.\qed

We have shown that the existence of DQE-metrics is necessary for the convergence of combinatorial CR-curvature flow. Actually, it is almost sufficient. In fact, we have
\theorem \label{main-proof}Suppose there is a DQE-metric $r^\ast$ such that $K^\ast=\lambda^* r^\ast$. Then the solution of the normalized CR-curvature flow $\dot{u}=L^T(K^\ast-K)$ exists for all time and converges to $r^\ast$ if the initial metric $r(0)$ is close to $r^\ast$, or equivalently, the initial normalized quadratic energy $\mathcal{C}(r(0))=\|K(0)-K^\ast\|^2=\sum_{i=1}^N (K_i(0)-K_i^\ast)^2$ is small enough.\\
\proof Denote $\Gamma(u)=L^T(K^\ast-K)$. The above normalized combinatorial CR-curvature flow is an autonomous ODE system $\dot{u}=\Gamma(u)$. We want to prove $u^\ast$ is a local attractor of the flow. By calculation we find $$\Gamma(u^\ast)=0,$$ and $$D_{u^\ast}\Gamma(u)=-L^TL \leq 0.$$ Although $-L^TL$ is only negative semi-definite in the whole space $\mathfrak{M}_{r}$, it is negative definite along the curvature flow due to the transversality announced in Remark \ref{transverse}. Therefore, $u^\ast$ is a local attractor of the flow. The system is asymptotically stable at $u^\ast$. The following four conditions
\begin{itemize}
  \item The initial metric $r(0)$ is close to $r^\ast$
  \item The initial metric $u(0)$ is close to $u^\ast$
  \item The initial CR-curvature $K(0)$ is close to $K^\ast$
  \item The initial normalized quadratic energy metric $\mathcal{C}(0)$ is small enough
\end{itemize}
are all equivalent due to the locally diffeomorphism between $K$ and $u$. They imply both the long time existence and convergence of combinatorial CR-curvature flow (\ref{Combina-Calabi-F}).\qed

Similarly, we can get
\corollary Given an admissible CR-curvature $\bar{K}$, that is, there is a metric $\bar{u}\in\mathfrak{M}_u$ that determines $\bar{K}=K(\bar{u})$. Then the solution of the modified CR-curvature flow $\dot{u}=L^T(\bar{K}-K)$ exists for all time and converges to $\bar{u}$ if the initial metric $r(0)$ is close to $\bar{r}$, or equivalently, the initial modified quadratic energy $\mathcal{C}(r(0))=\|K(0)-\bar{K}\|^2=\sum_{i=1}^N (K_i(0)-\bar{K_i})^2$ is small enough.\\

\remark Remember $\Lambda=LR^{-1}=\frac{\partial(K_{1},\cdots,K_{N})}{\partial(r_{1},\cdots,r_{N})}$ is positive semi-definite. Consider combinatorial curvature flow $\dot{r}=-\frac{1}{2}\nabla_{r} \mathcal{C}=-\Lambda K$. This flow is also a fourth order combinatorial CR-curvature flow. Along this flow, $\|r\|^2$ is invariant. Moreover, this flow has similar convergence properties with combinatorial CR-curvature flow (\ref{Combina-Calabi-F}).

\section{Combinatorial CR-curvature flow of second order}\label{2th}
\definition The second order combinatorial CR-curvature flow is
\begin{equation}
\dot{r}=-K.
\end{equation}
Consider the normalized flow
\begin{equation} \label{normalize-2th-flow}
\dot{r}=\lambda r-K,
\end{equation}
where $\lambda=\frac{S}{\|r\|^2}$. Since $\frac{d\|r\|^{2}}{dt}=2r^T\dot{r}=2r^T(\lambda r-K)=0$, then $\|r\|^{2}$ is a constant along the normalized flow (\ref{normalize-2th-flow}). By calculation, we get $\dot{K}=-\Lambda K$. $\Lambda$ is a second order discrete differential operator. That's why we call it flow of second order. Moreover, we have $\mathcal{\dot{C}}=-2K^T\Lambda K$, $\dot{S}=-\lambda S-\mathcal{C}=-\|K-\lambda r\|^2\leq 0$.

\theorem If the solution of combinatorial CR-curvature flow (\ref{normalize-2th-flow}) exists for all time and converges to a non-degenerate sphere packing metric $r(+\infty)$, then DQE-metric exists. Moreover, $r(+\infty)$ is indeed one.\\
\proof Along the normalized flow (\ref{normalize-2th-flow}), $\dot{S}=-\lambda S-\mathcal{C}=-\|K-\lambda r\|^2\leq 0$. Under the hypothesis that $r(t)\rightarrow r(+\infty)$ we know $\dot{S}(+\infty)=0$, which implies $K(+\infty)=\lambda r(+\infty)$.\qed

Let $\Gamma(r)=\lambda r-K=\frac{S}{\|r\|^2}r-K$, we have
\begin{equation*}
\begin{aligned}
D_{r}\Gamma(r)&=-D_{r}K+\lambda D_{r}r+r(D_{r}\lambda)^T\\
&=-\Lambda+\lambda I_{N}+r(\frac{\Lambda r+K}{\|r\|^{2}}-\frac{2S}{\|r\|^{4}}r)^{T}\\
&=\lambda I_{N}-\Lambda+\frac{rK^{T}}{\|r\|^{2}}-2S\frac{rr^{T}}{\|r\|^{4}}\\
&=\lambda (I_{N}-\frac{rr^{T}}{\|r\|^{2}})-\Lambda+\frac{r(K-\lambda r)^{T}}{\|r\|^{2}}.
\end{aligned}
\end{equation*}
At DQE-metric point $r^*$ we have
\begin{equation}
D_{r}\Gamma\big|_{r^*}=\lambda^* (I_{N}-\frac{rr^{T}}{\|r\|^{2}})-\Lambda.
\end{equation}

Next we give two sufficient conditions so as the normalized flow (\ref{normalize-2th-flow}) exists for all time and converges to DQE-metric.

\theorem \label{nonpositive-condition}
Given $(M^3,\mathcal{T})$. Let $r^*$ be a DQE-metric with $K^*=\lambda^* r^*, \lambda^*\leq 0$. Assuming the initial normalized quadratic energy $\mathcal{C}(r(0))=\|K(0)-K^\ast\|^2$ is small enough, the solution of the normalized flow (\ref{normalize-2th-flow}) exists for $t\in [0,+\infty)$ and converges to $r^*$. \\
\proof Notice that $I_{N}-rr^T/\|r\|^{2}$ is always positive semi-definite. As $\lambda^*\leq 0$ at $r^*$,
$$D_{r}\Gamma\big|_{r^*}=\lambda^* (I_{N}-\frac{rr^{T}}{\|r\|^{2}})-\Lambda\leq0.$$
Moreover, $rank\,(D_{r}\Gamma\big|_{r^*})=N-1$ and the kernel of $D_{r}\Gamma\big|_{r^*}$ is exactly $tr,\,t\in \mathds{R}$. Along the normalized flow (\ref{normalize-2th-flow}), $\|r\|$ is a constant and never scales. Think $D_{r}\Gamma\big|_{r^*}$ as a negative positive definite matrix, hence each flat or negative DQE-metric $r^*$ is a local attractor of the normalized flow (\ref{normalize-2th-flow}).\qed

\theorem
Given $(M^3,\mathcal{T})$. Let $r^*$ be a DQE-metric with $K^*=\lambda^*r^*$. Assuming $\lambda_1(\Lambda)$, the first eigenvalue of $\Lambda$ at $r^*$, satisfies
$\lambda_1(\Lambda)>\lambda^*$. Furthermore, assuming the initial normalized quadratic energy $\mathcal{C}(r(0))=\|K(0)-K^\ast\|^2$ is small enough, the solution of the normalized flow (\ref{normalize-2th-flow}) exists for $t\in [0,+\infty)$ and converges to $r^*$. \\
\proof Select an orthogonal matrix $P$, such that $P^T\Lambda P=diag\{0,\lambda_1,\cdots,\lambda_{N-1}\}$. Write $P=(e_0,e_1,\cdots,e_{N-1})$,
where $e_i$ is the $(i+1)$-column of $P$. Then $\Lambda e_0=0$ and $\Lambda e_i=\lambda_i e_i,\,1\leq i\leq N-1$,
which implies $e_0=r/\|r\|$ and $r\perp e_i,\,1\leq i\leq N-1$.
Then $\big(I_{N}-\frac{rr^T}{\|r\|^{2}}\big)e_i=e_i$, $1\leq i\leq N-1$, which implies $P^T\big(I_{N}-\frac{rr^T}{\|r\|^{2}}\big)P=diag\{0,1,\cdots,1\}$. Therefore,
\begin{equation*}
D_{r}\Gamma\big|_{r^*}=\lambda^* (I_{N}-\frac{rr^{T}}{\|r\|^{2}})-\Lambda=P^Tdiag\{0,\lambda^*-\lambda_1,\cdots,\lambda^*-\lambda_{N-1}\}P.
\end{equation*}
If $\lambda_1(\Lambda)>\lambda^*$, then $D_{r}\Gamma\big|_{r^*}\leq0$, $rank\,(D_{r}\Gamma\big|_{r^*})=N-1$. The proof can be finished completely by the same reason explained in the proof of Theorem \ref{nonpositive-condition}.
\qed

\section{Combinatorial G-Laplacian and G-curvature flow}
\label{G-flow}
Besides CR-curvature $K_i$, there are other types of combinatorial scalar curvatures. In \cite{G5}, Glickenstein considered the curvature $C_i=K_ir_i$. From now on, we call $C_i$ ``G-curvature" for short. CR-curvature is invariant under metric scaling, while G-curvature scales in the appropriate way. G-curvature is like scalar curvature times the volume measure, so scales like length.

The total curvature functional becomes $S=\sum_{i=1}^N K_{i}r_{i}=\sum_{i=1}^N C_i.$ Define a new discrete quadratic energy $$\mathcal{E}=\|C\|^2=\sum_{i=1}^NC_i^2.$$ Now we consider its negative gradient flow
\begin{equation}
\dot{u}(t)=-\frac{1}{2}\nabla_{u} \mathcal{E}=-G^TC
\end{equation}
Here $C=(C_1,\cdots,C_N)^T$, $G=\frac{\partial(C_{1},\cdots,C_{N})}{\partial(u_{1},\cdots,u_{N})}$, $\nabla_{u} \mathcal{E}=2G^TC$.
Therefore
$$\nabla_{u} S=C,$$
$$Hess_{u} S=G.$$
The matrix $G$ is symmetric since it is the Hessian of the total scalar curvature functional (or see \cite{G5}). We call it G-Laplacian. Remember that $R=diag\{r_1,\cdots,r_N\}$, $\Lambda=\frac{\partial(K_{1},\cdots,K_{N})}{\partial(r_{1},\cdots,r_{N})}$. Denote $\Omega=diag\{C_1,\cdots,C_N\}$, then by calculation we obtain $$G=\Omega+\widetilde{L}=\Omega+R\Lambda R.$$
\lemma
If $C_i\geq0$ for any $i$ and there is at least one $C_i>0$, then the matrix $G$ is positive definite.\\
\proof
The matrix $\Omega$ and $R\Lambda R$ are both semi-positive definite. The kernel of $R\Lambda R$ is constant column vector $\textbf{1}=(1,\cdots,1)^T$ up to scaling, however, this vector makes $\textbf{1}^T\Omega\textbf{1}>0$.\qed

When $G$ is positive definite, good properties will come out locally. Let $$\mathfrak{M}_{r}^\ast=\left\{\;r\in\mathfrak{M}_{r}\;\big|\;G(r)>0\right\}.$$
All metrics in this open set determine positive definite G-Laplacians exactly. We use $\mathfrak{M}_{u}^\ast$ to denote its corresponding domain under coordinate transformations $u_i=\ln r_i$.

For any convex set $U\subset \mathfrak{M}_{u}^\ast$, since $G=Hess_{u} S$ is positive definite, the G-curvature map $C(u)=\nabla_{u} S$ is an embedding on $U$. This means locally every G-curvature determines a unique metric when G-curvature is positive definite. G-curvature contains the $r$ term, hence is not invariant under scaling. This local rigidity property of G-curvature takes scaling into account, and is different with the case of Cooper and Rivin's curvature.

Now we give an example which is an analogue of shrinking sphere in the case of smooth Ricci flow.
\example
Given $(M^3,\mathcal {T})$, and $r(0)$ is the initial sphere packing metric. Let $C(0)$ be the corresponding initial G-curvature and $G(0)$ the initial G-Laplacian. If $G(0)C(0)$ is a constant column vector, then the combinatorial G-curvature flow has a shrinking or steady solution. Suppose $G(0)C(0)=\lambda \textbf{1}$, since $\lambda N=\textbf{1}^TG(0)C(0)=\textbf{1}^T\Omega(0)C(0)=C(0)^TC(0)=\mathcal{E}(0)$, we get
$$\lambda=\frac{\mathcal{E}(0)}{N}\geq0.$$ Then the G-curvature flow $\dot{u}=-G^TC$ with initial metric $u(0)$ has a unique solution
$$r(t)=\frac{r(0)}{\sqrt{1+2\lambda t}},\;\;-\frac{1}{2\lambda}<t<+\infty.$$

In fact we can solve the above G-curvature flow equation by method of undetermined coefficients. The solution is shrinking when $\lambda>0$ and steady when $\lambda=0$. Notice that, the condition that $G(0)C(0)$ is a constant column vector can be interpreted as ``isotropy". It is an analogue of the standard ``round" sphere in a certain sense. \qed

It is conceivable the G-curvature flow will exhibit good local convergence properties when the G-Laplacian is positive definite. For any user prescribed target G-curvature with positive G-Laplacians, the modified G-curvature flow tends to find the target G-curvature automatically.
\theorem
Given an admissible G-curvature $\bar{C}$ with positive G-Laplacian, that is, there is a metric $\bar{u}\in \mathfrak{M}_{u}^\ast$ that determines $\bar{C}=C(\bar{u})$. Consider the modified quadratic energy $$\mathcal{E}=\sum_{i=1}^N(C_i-\bar{C_i})^2.$$ The modified combinatorial G-curvature flow
\begin{equation}
\dot{u}(t)=-\frac{1}{2}\nabla_{u} \mathcal{E}=G^T(\bar{C}-C)
\end{equation}
exists for all $t\geq0$ and converges to $\bar{u}$ as $t\rightarrow +\infty$ whenever the initial quadratic energy is small enough.\\
\proof
Notice that the G-Laplacian $\bar{G}=G(\bar{u})$ is positive definite, then use the methods in the proof of Theorem \ref{main-proof} we can get the above results. \qed
\remark For any user prescribed target G-curvature $\bar{C}$ (suppose $\bar{C}$ is well selected, \emph{i.e.} it is admissible and it determines positive definite G-Laplacian), any algorithm aiming to minimize $\mathcal{E}=\sum_{i=1}^N(C_i-\bar{C_i})^2$, such as Newton's gradient descent algorithm, tends to find sphere packing metric $\bar{r}$ automatically.
\section{Hyperbolic geometric background}
For hyperbolic manifolds, it's more natural to built the triangulation by hyperbolic tetrahedrons. Each block is embedded in $\mathds{H}^3$. For a hyperbolic conformal tetrahedron $\{i,j,k,l\}\in T$, Cooper-Rivin showed that the matrix $$\Lambda_{ijkl}=-\frac{\partial(\alpha_{i},\alpha_{j},\alpha_{k},\alpha_{l})}{\partial(r_{i},r_{j},r_{k},r_{l})}$$ is positive definite in \cite{CR} and \cite{Ri}. So $$\Lambda=\sum_{\{i,j,k,l\}\in T}\Lambda_{ijkl}$$ is also positive definite due to the same reason in the proof of Lemma \ref{semi-positive}. Notice that the total curvature functional is changed to $$S=2vol(M)+\sum_{i=1}^NK_ir_i.$$ Using Schl$\ddot{a}$fli formula again, we have $$dS=\sum_{i=1}^NK_idr_i+\left(2dvol(M)+\sum_{i=1}^Nr_idK_i\right)=\sum_{i=1}^NK_idr_i,$$ and hence $$\nabla_{r}S=K, \;Hess_{r}S=\Lambda.$$ Therefore, the curvature map $K=K(r)$ is also an embedding and hence locally a diffeomorphism. This implies the following
\theorem In hyperbolic geometric background, the sphere packing metric which determines constant zero curvature is isolated.\qed

Consider the following three discrete curvature flows
\begin{itemize}
  \item $\dot{r}=-K$,
  \item $\dot{r}=-\frac{1}{2}\nabla_{r} \mathcal{C}=-\Lambda K$,
  \item $\dot{u}=-\frac{1}{2}\nabla_{u} \mathcal{C}=-L^T K$.
\end{itemize}
Here $L=\Lambda R_h$, $R_h=diag\left\{ \sinh(r_1),\cdots,\sinh(r_N)\right\}$ and the coordinate change is $u_i=\ln tanh \frac{r_i}{2}$. It's interesting that these three discrete curvature flows possess similar convergence properties. We only list the results here without proof.

\theorem If the solution of any of the above three curvature flows exists for all time and converges to a non-degenerate metric $r(+\infty)$, then there exists constant zero curvature metric and $r(+\infty)$ is indeed one.
\theorem Suppose there is a constant zero metric $r^\ast$, then these three discrete curvature flows exist for all time and converge to $r^\ast$ when the initial energy is sufficiently small.

\remark
As was discussed in section \ref{G-flow}, we can also investigate G-curvature flow in hyperbolic geometric background. In this case, the G-curvature should be $C_i=K_i \sinh r_i$.
We omit the details since the processes and the conclusions are almost the same as before.

\section{Future Work}
There are something more we want to do in the future.
\subsection{Find combinatorial condition so as DQE-metric exists}
We do not even know whether DQE-metrics exist for a triangulated compact manifold $M^3$. Inspired by Thurston's work in \cite{T1} (or see \cite{CL1}), there are topological and combinatorial obstruction for the existence of constant curvature metric. We want to find similar topological and combinatorial condition for the existence of DQE-metrics. We also want to know how many DQE-metrics there are in a conformal class.
\subsection{How to define discrete curvatures in high dimensions? Dose discrete curvatures carry any topological information?}
Consider a $n$-dimensional triangulated manifold $(M,\mathcal{T})$ with a given flat cone metric, such that it becomes a PL-manifold. If $n=2$, the only singularity arises at vertices. However, in the case $n\geq 3$, singularities arise both at vertices, edges and even $(n-2)$-dimensional faces. What is the ``right" analog of the sectional curvature, Ricci curvature and scalar curvature? It seems the ``right" definition of curvatures is variant according to different purpose, such as classification purpose, capturing topological information purpose, finding good metric purpose and so on. There are not solely one ``right" definition of curvature.
\subsection{How to estimate $\lambda_1(\Lambda)$?}
We proved that the DQE-metrics $r^*$ with $\lambda_1(\Lambda)>\lambda^*$ is a local attractor of the second order combinatorial CR-curvature flow (\ref{normalize-2th-flow}). We want to develop efficient methods to estimate the lower bound of $\lambda_1(\Lambda)$. We also want to know if there are topological and combinatorial obstructions for $\lambda_1(\Lambda)>\lambda^*$.
\subsection{How far can the ``flow" method go?}
We still consider a $n$-dimensional triangulated manifold $(M,\mathcal{T})$. Suppose we have defined a type of scalar curvature $K: V\rightarrow \mathds{R}$ for each sphere packing metric $r$. Then $K$ determines a discrete energy $\mathcal{C}(r)=\sum_{i=1}^N K_i^2$. Similarly, we consider the negative gradient flow of this energy
$$\dot{u}=-\frac{1}{2}\nabla_u \mathcal{C}=-L^TK.$$
Then $$\dot{\mathcal{C}}=-2K^TLL^TK\leq 0.$$
Notice that $L=\frac{\partial(K_{1},\cdots,K_{N})}{\partial(u_{1},\cdots,u_{N})}$ and $rank(LL^T)=rank(L)$. If the combinatorial scalar curvature $K$ defined to be have the property $rank(L)=N-1$, then this flow converges to $Ker(L^T)$ when the initial energy is sufficiently close to $Ker(L^T)$.

The properties of the matrix $L$ is of most importance. The matrix $L$ came from the relation between curvature $K$ and the metric $r$ (or $u$). The full rank (up to scaling) of $L$ is equivalent to locally rigidity, $i.e.$ curvature determines metric. If a curvature $K$ is well defined and can be used for classifying PL-manifolds, then it should determine the PL-metric in a certain sense. As explained above, ``classification", ``rigidity", ``full rank of $L$ up to scaling" are closely related.

\subsection{Combinatorial Gauss-Bonnet formula in higher dimensions}
Suppose the discrete scalar curvatures defined in a certain manner have the properties that they are homogeneous functions of degree zero. Then they are invariant under scaling of metrics. Three are $N$ vertices and $N$ curvatures $K_1,...,K_N$, however, the essentially independent number of metrics $r_1,...,r_N$ is $N-1$ when modular scaling. Thus the curvatures $K_1,...,K_N$ must have some relations, that is, $K$ lies in a no more than $(N-1)$-dimensional submanifold in $\mathds{R}^N$. If the curvatures possess locally rigidity additionally, they must lie in a hypersurface in $\mathds{R}^N$. For the case $dim(M)=2$, people have already known the combinatorial Gauss-Bonnet formula $$\sum_{i=1}^NK_i=2\pi \chi(M).$$ Is there an analogue of combinatorial Gauss-Bonnet formula in higher dimensions? Can we write the hypersurface out clearly?

\subsection{How to get global rigidity?}
The space of admissible sphere packing metrics $\mathfrak{M}_{r}$ is not convex, which is the main obstruction to get global rigidity. In \cite{L2}, Feng Luo suggested finding a good parametrization so as the admissible space of metrics became convex in the new coordinate. In \cite{CR}, Cooper-Rivin considered a new coordinate system $k_i=1/r_i$. The admissible space of metrics became convex under this coordinate change. Unfortunately, negativeness of the Hessian matrix of total curvature functional under new coordinate seems hard to preserve. It seems not easy to find a better parametrization so that the admissible space of metrics becomes convex and the positiveness (negativeness) of the Hessian matrix is preserved simultaneously. We still need to think up more efficient methods to get global properties.\\[8pt]

{\Acknowledgment Acknowledgments}\\[8pt]
The first author would like to show his greatest respect to Professor Gang Tian who brought the author to the area of combinatorial curvature flows. The second author would like to thank Professor Zhang Xiao for the invitation to AMSS. Both authors would like to give special thanks to Professor Guanxiang Wang for reading the paper carefully and giving numerous improvements to the paper. Both authors would also like to thank Professor Glickenstein, Feng Luo, Yuguang Shi and Meiyue Jiang for many helpful conversations. Both authors would like to give special thanks to Dr. Zhou Da, Wang Ding for their encouragements.

(Huabin Ge) Beijing International Center for Mathematical Research, Peking Univ., Beijing 100871, P.R. China

E-mail: gehuabin@pku.edu.cn\\[2pt]

(Xu Xu) School of Mathematics and Statistics, Wuhan Univ., Wuhan 430072, P.R. China

E-mail: xuxu2@whu.edu.cn\\[2pt]

\end{document}